\documentclass{article}

\usepackage{latexsym}
\usepackage{graphicx}
\usepackage{bm}% bold math

\begin{document}

 \def\prob#1{\mbox{\,{\bf Pr}$[#1]$}}

\title{PROPAGATION TIME IN STOCHASTIC COMMUNICATION NETWORKS
}
\author{JONATHAN E. ROWE \\ School of Computer Science \\ University of Birmingham \\ Birmingham B15 2TT \\ J.E.Rowe@cs.bham.ac.uk \\ \\
BORIS MITAVSKIY \\ Division of Genomic Medicine \\
School of Medicine \\
University of Sheffield \\
Sheffield S10 2JF \\ B.Mitavskiy@sheffield.ac.uk}

\maketitle

\begin{abstract}
Dynamical processes taking place on networks have received much attention in recent years, especially on various models of random graphs (including Òsmall worldÓ and Òscale freeÓ networks). They model a variety of phenomena, including the spread of information on the Internet; the outbreak of epidemics in a spatially structured population; and communication between randomly dispersed processors in an ad hoc wireless network. Typically, research has concentrated on the existence and size of a large connected component (representing, say, the size of the epidemic) in a percolation model, or uses differential equations to study the dynamics using a mean-field approximation in an infinite graph. Here we investigate the time taken for information to propagate from a single source through a finite network, as a function of the number of nodes and the network topology. We assume that time is discrete, and that nodes attempt to transmit to their neighbours in parallel, with a given probability of success. We solve this problem exactly for several specific topologies, and use a large-deviation theorem to derive general asymptotic bounds, which apply to any family of networks where the diameter grows at least logarithmically in the number of nodes. We use these bounds, for example, to show that a scale-free network has propagation time logarithmic in the number of nodes, and inversely proportional to the transmission probability.
\end{abstract}

%\keywords{network dynamics; propagation time; network topology}

%\pacs{02.50.Fz,89.75.Hc,05.45.-a,87.23.-n}

%{\bf Summary} We investigate the time taken for information to propagate through a network, as a function of the network size and topology, assuming that copying occurs in discrete time steps, in parallel, and with limited reliability.

\section{Introduction}

Within a few years we will be able to produce vast numbers of microscopic, extremely cheap, computer processors~\cite{abelson:2000}. These could be randomly distributed (or painted) on a surface and, by making use of their massive parallelism, form an intelligent, computational lawn. However, each processor will only be able to communicate over a short range and with limited reliability. An obvious question is: how long would a message take to spread across the network starting from a single source? Similar questions arise in epidemiology~\cite{keeling:2005, barthelemy:2004, boguna:2003, grassberger:1983, newman:2002, sander:2002}. Given a spatially distributed population, in which individuals infect their neighbours with a certain probability, how long before the whole population is infected? 

\begin{figure}
\begin{center}
\includegraphics[scale=0.5]{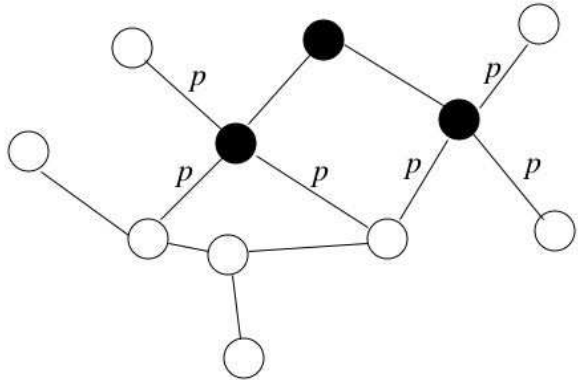}
\end{center}
\caption{Nodes that have the information (black) try, in parallel, to copy to their neighbours. Transmissions succeed with probability $p$.}
\label{fig1}
\end{figure}

There are several different questions that can be posed about propagation in networks. If there is only one chance for a node to successfuly transmit to its neighbour, one can ask under what conditions (and with what probability) most or all of the network is ``infected'' in the steady state. This requires the investigation of the percolation structure of the system~\cite{grassberger:1983,sander:2002,benjamini:2003,pemantle:1995,pemantle:1994,grimmett:1989}. Alternatively, one can focus on the dynamics. A typical approach here is to assume an infinite network with sufficient symmetry to make use use of a mean-field approximation in continuous time~\cite{keeling:2005,barthelemy:2004,newman:2002}. In this paper, we wish to investigate propagation through a \emph{finite} network of arbitrary topology (including so-called \emph{scale-free} networks such as the Internet~\cite{watts:2004,song:2005,barabasi:2001}), assuming that transmissions take place in parallel, in discrete time steps, and with some fixed probability of success. Transmissions are attempted at every time step until successful.

\begin{figure}
\begin{center}
\includegraphics[scale=0.45]{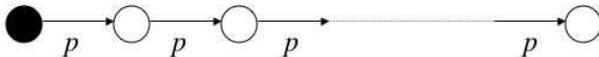}
\end{center}
\caption{A simple chain of $n$ nodes.}
\label{fig2}
\end{figure}

Suppose we have a network of computers and some information is located on one of them. At each time step, processors with the information try to copy it to their neighbours, with success probability $p$ (see figure~\ref{fig1}). We wish to estimate the expected time for the information to spread to all nodes, which we call the \emph{propagation time}, $E[n]$, on a network containing $n$ nodes. Similarly, one could consider the spread of an infectious disease~\cite{anderson:1991}, where $p$ is the infection rate, or the spread of a mutant gene in a metapopulation~\cite{lieberman:2005,hanski:2004}. 
The simplest network to consider is a chain of $n$ nodes (see figure~\ref{fig2}). If the information starts at one end, then the expected time, $E[n]$, for it to have crossed the network is $(n-1)/p$. This result can be derived from a recurrence relation for the propagation time given $k$ remaining nodes: 
\begin{equation}
	E[k] = 1 + pE[k-1] + (1-p)E[k]
\end{equation}
A similar recurrence enables us to solve the case of a ring of $n$ nodes. The exact result is complicated but to a good approximation is $(n-1)/2p$.

\section{The general recurrence equation}

Such a recurrence can be derived for any network: the propagation time starting from a particular situation can be broken down into the possible cases occurring after a single time step, with their associated probabilities.

We consider the general problem of a random sequence $X_1, X_2, \ldots$ of states from some set $A$, and a subset $D \subseteq A$ of desired states. We can derive a recurrence relation for the first hitting time of the desired subset in terms of all the situations that could possibly arise after a single time step, and their associated probabilities. The first hitting time is defined as $T = \min \{t \, | \,X_t \in D\}$. In our case, the random sequence comes from the different states of the network as the information (or infection) is propagated from node to node. The desired state is when all the nodes have been infected.

Suppose after the first time step the network is in state $X_1 = k$.  We consider the \emph{conditional} probability space that arises from this situation. We write $E_k$ to denote expectation with
respect to this conditional probability space, and consider the shifted
stochastic process 
\[
	Y_1 = X_2, Y_2 = X_3, \ldots, Y_n = X_{n+1}, \ldots
\]  
Let $T_k = \min \{t \, | \, Y_t \in D\}$ be the first hitting time after this first step has been made.
We then have: 
\[
E[T] = \sum_{k \in A} \prob{X_1 = k} \cdot E_k[T_k] + 1
\]

{\bf Proof}
\begin{eqnarray*}
  E[T] & = &  \sum_{j = 1}^{\infty}  \prob{T = j} \cdot j \\
             & = & \sum_{j = 1}^{\infty}  \left(\sum_{k \in A} \prob{X_1 = k} \cdot \prob{T = j \, | \, X_1 = k}\right) \cdot ((j-1)+1) \\
             & = & \sum_{k \in A} \prob{X_1 = k} \left(\sum_{j = 1}^{\infty} \prob{T = j \, | \, X_1 = k} \cdot (j-1)+1\right) \\
             & = & \sum_{k \in A} \prob{X_1 = k} \left(\sum_{j =1}^{\infty} \prob{T = j \, | \, X_1 = k} \cdot (j-1)\right)+\sum_{k \in A} 
             \prob{X_1 = k} \\
             & = & \sum_{k \in A} \prob{X_1 = k} \cdot E_k[T_k] + 1
\end{eqnarray*}
\hfill
$\Box$

\section{The hub model}

A more complex example is the \emph{hub} in which the information starts at a central node which repeatedly tries to transmit to $n$ neighbours. For example, consider a transmitter signalling to a number of receivers, or a collection of $n$ people that independently have probability $p$ of contracting a disease. The recurrence relation becomes: 
\[	 E[n] = 1 + \sum_{k=0}^{n} {n \choose k} p^{n-k} q^k E[k]
\]
where $q = 1-p$. We rearrange to give the recurrence relation:
\[
(1 - q^n) E[n] = 1 + \sum_{k=0}^{n-1} {n \choose k} p^{n-k} q^k E[k]
\]
with $E[0] = 0$. The transmission probability is $p$ and $q = 1-p$.  We claim that $E[n] = \Theta(\log(n+1))$.

We first prove the general result that, if $k$ is distributed according to the Binomial distribution $B(n, q)$, then the expected value of $\log(k+1)$ is $\Theta(\log(n+1))$. 
%In particular,
%\[
%	\log(n+1) - 1/q \leq \sum_{k=0}^n {n \choose k} p^{n-k} q^k \log(k+1) \leq \log(n+1) - \log\left(\frac{2}{1+q}\right)
%\]
%(assuming natural logs --- for other bases, the constant in the lower bound has to change accordingly). 

{\bf Proof} 
To prove the upper bound, first note that, for $0 < q < 1$,  and for all $n \geq 1$ 
\[
  nq + 1 \leq \frac{(n+1)(1+q)}{2}
\]
Then, by concave property of logs,
\begin{eqnarray*}
 & & \sum_{k=0}^n {n \choose k} p^{n-k} q^k \log(k+1)  \\ & \leq & \log\left( \sum_{k=0}^n {n \choose k} p^{n-k} q^k (k+1) \right) \\
& = & \log(nq+1) \\
& \leq & \log\left( \frac{(n+1)(1+q)}{2} \right) \\
& = & \log(n+1) - \log\left(\frac{2}{1+q}\right)
\end{eqnarray*}

Now note that the lower bound holds, if and only if:
\[
  \sum_{k=0}^n {n \choose k} p^{n-k} q^k \log \left(\frac{n+1}{k+1}\right) \leq 1/q
\]
which we prove as follows:
\begin{eqnarray*}
 & &  \sum_{k=0}^n {n \choose k} p^{n-k} q^k \log \left(\frac{n+1}{k+1}\right) \\ & = &    \sum_{k=0}^n {n \choose k} p^{n-k} q^k \log \left(1 + \frac{n-k}{k+1}\right) \\
  & \leq &   \sum_{k=0}^n {n \choose k} p^{n-k} q^k \left(\frac{n-k}{k+1}\right) \\
  & = &   \sum_{k=0}^n {n \choose k} p^{n-k} q^k \left(\frac{n}{k+1}\right) -   \sum_{k=0}^n {n \choose k} p^{n-k} q^k \left(\frac{k}{k+1}\right) \\
  & \leq & n \sum_{k=0}^n {n \choose k} p^{n-k} q^k \left(\frac{1}{k+1}\right) \\ 
  & = & \frac{n}{n+1} \sum_{k=0}^n  {n+1 \choose k+1} p^{n-k} q^k \\
  & = & \frac{n}{n+1} \sum_{k = 1}^{n+1} {n+1 \choose k} p^{n-k+1} q^{k-1} \\
  & \leq & \frac{1}{q} \cdot \frac{n}{n+1} \sum_{k = 0}^{n+1} {n+1 \choose k} p^{n+1-k} q^{k} \\
  & \leq & \frac{1}{q}
\end{eqnarray*}
(where we used the fact that $\log(1+x) \leq x$ for all $x \geq 0$).
\hfill
$\Box$

We now use this result to show that
the propagation time for a hub with $n$ clients is $E[n] = \Theta(\log(n+1))$. In particular,
\[
    q \log(n+1) \leq E[n] \leq \left( \frac{1}{\log\frac{2}{1+q}} \right) \log(n+1)
\]

{\bf Proof} 
We prove, by induction, that $E[n] \leq A \log(n+1)$, where
\[
  A = \frac{1}{\log\frac{2}{1+q}}
\]

The case $n=0$ is easy, since $E[0] = 0 = A \log 1$. Now suppose $n \geq 1$ and that the hypothesis is true for all $0 \leq k < n$. Then
\begin{eqnarray*}
 & & (1-q^n) E[n]  \\ & \leq 1 & + \sum_{k=0}^{n-1} {n \choose k} p^{n-k} q^k A \log(k+1)\\
& = & 1 - A q^n \log(n+1) + A \sum_{k=0}^{n} {n \choose k} p^{n-k} q^k  \log(k+1)\\
& \leq & 1 - A q^n \log(n+1) + A\log(n+1) - A \log \left(\frac{2}{1+q}\right) \\
& = & A (1 - q^n) \log (n+1)
\end{eqnarray*}
and the result follows.
 
Secondly, we show  by induction that $E[n] \geq q \log(n+1)$. The case $n=0$ is again easy. Now suppose the hypothesis is true for all $0 \leq k <  n$. Then
\begin{eqnarray*}
 & & (1-q^n) E[n] \\ & \geq & 1 + q \sum_{k=0}^{n-1} p^{n-k} q^k \log(k+1) \\
& = & 1 - q^{n+1} \log(n+1) + q \sum_{k=0}^{n} p^{n-k} q^k \log(k+1) \\
& \geq & 1 - q^{n+1} \log(n+1) + q\left( \log(n+1) - \frac{1}{q}\right) \\
& = & q(1 - q^n)\log(n+1)
\end{eqnarray*}
%and the result follows.
\hfill
$\Box$

\section{Epidemiology and perfect mixing}

 In epidemiology models of the spread of infectious diseases, it is common to assume \emph{perfect mixing}: that every individual interacts with every other individual. This corresponds to having a complete graph, with every node connected to every other node. The propagation time for a complete graph is bounded by a constant: it does not depend on the number of nodes. Moreover, as $n$ gets large, the propagation time is just two time steps (with increasingly high probability). 
Suppose we have a complete graph on $n+1$ nodes and initially one node is infected. The probability that the infection passes to a neighbour is $p$ and we let $q = 1-p$. After a single time step it is very likely that close to $n p$ new nodes have been infected. In fact, Chernoff's inequality~\cite{mitzenmacher:2005} tells us that the probability that less than $(1 - \delta) n p$ nodes have been infected is less than $\exp(- np\delta^2/2)$, where we can make $\delta$ as small as we like. This means that, the more nodes there are, the surer we can be that nearly $np$ nodes are infected after one time step. The probability that the remaining nodes get infected on the next time step is therefore close to $(1 - q^{np})^{nq}$. Using the fact that $(1-x)^n \geq 1-nx$ for all $0\leq x\leq 1$, we see that
$(1 - q^{np})^{nq} \geq 1 - nq^{np + 1} \geq 1 - r^n 
$
(for some $r$ in the range $0 < q^{np} < r < 1$), which is $1 - e^{-O(n)}$.

Similarly for a complete bipartite graph, with $n$ nodes in each set, the probability that a single node in one set infects $np$ nodes in the second gets arbitrarily close to 1 as $n$ gets large. This is then enough to infect all the other nodes in the first set (again with arbitrarily high probability). Then, on the third time step, the remaining nodes of the second step get infected. This analysis can be extended to complete multipartite graphs in an obvious way.

In our model, nodes are either infected or not. This corresponds to the SI model of epidemiology (Susceptible-Infected~\cite{anderson:1991}). If $I$ is the number of infected people and $S=n=I$ are the remaining susceptibles, we would like to know how many more people become infected in one time step. For the complete graph (perfect mixing) the number of newly infected people is binomially distributed between $0$ and $S$ with success probability $1 - (1-p)^I$. If $p$ is small, this is approximately equal to $pI$ and so the expected increase in infected people is close to $\Delta I \approx pSI$, which agrees with the standard SI model. To extend our model to more realistic scenarios, one would have to introduce a third state R (removed) for those people who can no longer be infected (due to immunity or death).

\begin{figure}
\begin{center}
\includegraphics[scale=0.3]{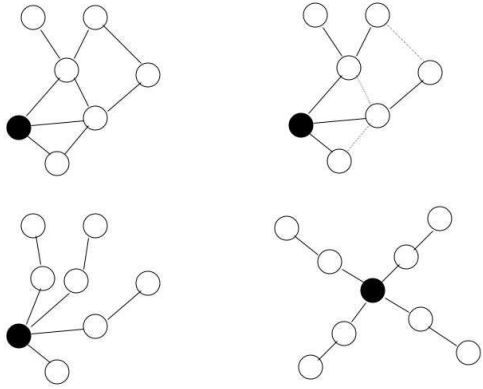}
\end{center}
\caption{To find an upper bound for the propagetion time of a network, we first find a minimum spanning tree. From this, we create a star graph, which we then balance.}
\label{fig3}
\end{figure}

\section{General upper and lower bounds}

For the general case, a lower bound on the propagation time, starting from a particular node, is given by the eccentricity of that node (divided by $p$). That is, the distance from the source to the most distant node in the network. This is because at least that number of successful transmissions will have to be made for the whole network to be infected. It is also possible to derive a general upper bound on the propagation time for an arbitrary network. The idea is to replace the network with a minimum spanning tree, rooted at the starting node. The propagation time on the tree is slower than for the original network, since we may have lost a number of ``short-cuts''. We then replace the tree with a star graph, with a hub at the starting node, and $b$ branches: one for every leaf of the tree. The length of each branch is the eccentricity, $d$, of the hub (see figure~\ref{fig3}). Using Chernoff bounds, we can prove  that the propagation time for a star graph is $O((d + \log b)/p)$.

To show this, we return to the general problem of estimating the first hitting time of a desired subset of states in a random sequence.
Let $A$ be a set of states and let $D \subseteq A$ be the desired subset. Let $Y(1), Y(2), \ldots$ be a random sequence of states from $A$ that satisfies the following \emph{monotonicity} properties:
\begin{enumerate}
\item If $Y(t) \in D$, then for all $k > 0$, $Y(t+k) \in D$.
\item $ \prob{Y(t+k) \in D | Y(t) \notin D} \geq \prob{Y(k) \in D}$
for all $t, k > 0$.
\end{enumerate}
In other words, the probability of reaching the desired state after a given time interval always improves, and once it is reached, it is never left. In the case of the star graph, the random variables $Y(t)$ will represent the minimum number of infected nodes along each branch at time $t$. The desired state is that all nodes on all branches are infected. Clearly, the probability of achieving this state in a fixed amount of time can only improve as time goes by, so the monotonicity condition is satisfied.

Define the first hitting time of the sequence to be $T(D) = \min \{t | Y(t) \in D\}$. 
Suppose we can find a time $\tau$ and constant $0 \leq \varepsilon < 1$ such that
\[
    \prob{Y(\tau) \in D} \geq 1 - \varepsilon
\]
Then we claim that
\[
    E[T] \leq  \frac{\tau}{1-\varepsilon} 
\]
That is, if after some time $\tau$ (which will in general depend on the structure of the problem), we have some lower bound on the success probability at that time, then we can use this fact to estimate the first hitting time for the whole process.

{\bf Proof}
According to the definition of expectation we have

\begin{eqnarray*}
% & & E[T]  =  
\sum_{k=1}^{\infty} k \prob{T=k} 
& = &  \left(\sum_{k=1}^{\tau} k \prob{T=k}\right) + \left( \sum_{k=\tau+1}^{2\tau} k \prob{T=k}\right) + \left(\sum_{k=2\tau+1}^{3\tau} k \prob{T=k}\right) + \cdots \\
& \leq & \left(\tau \sum_{k=1}^{\tau} \prob{T=k}\right) + \left(2\tau \sum_{k=\tau+1}^{2\tau} \prob{T=k} \right) +\left( 3\tau  \sum_{k=2\tau+1}^{3\tau} \prob{T=k}\right) + \cdots \\
& = & \tau \left(\sum_{k=1}^{\infty} \prob{T=k}\right) +  \tau \left(\sum_{k=\tau+1}^{\infty} \prob{T=k} \right) + \tau \left(\sum_{k=2\tau+1}^{\infty} \prob{T=k}\right) + \cdots \\
& = & \tau\left(1 + \prob{Y(\tau) \notin D} + \prob{Y(2\tau) \notin D} + \cdots \right)
\end{eqnarray*}

Now the second monotonicity condition can be equivalently stated as: 
\[
	\prob{Y(t+k) \notin D | Y(t) \notin D} \leq \prob{Y(k) \notin D}
\] 
for all $t, k > 0$. Using the definition of conditional probability, this gives us:
\[
	\prob{Y(t+k) \notin D} \leq \prob{Y(t) \notin D}\prob{Y(k) \notin D}
\]
We already know that $\prob{Y(\tau) \notin D}  \leq \varepsilon$. And by induction on $m$ we get  $\prob{Y(m\tau) \notin D}  \leq \varepsilon^m$. Therefore
\[
E[T] \leq \tau \sum_{m=0}^{\infty} \varepsilon^m = \frac{\tau}{1 - \varepsilon}
\]
as required.
\hfill
$\Box$

Now consider propagation in a star graph with $b$ branches of depth $d$. The problem is equivalent to $b$ parallel repeating Bernoulli trials $X_1(t), X_2(t), \ldots, X_b(t)$, each with success probability $p$. $X_j(t)$ is the number of infected nodes on branch $j$ at time $t$. We want the expected time until all of them have achieved at least $d$ successes. So let $Y(t) = \min_{1 \leq i \leq b} X_i (t)$, which certainly satisfies the monotonicity requirement.
Then
\begin{eqnarray*}
 & & \prob{Y(t)  <  (1-\delta)tp} \\ & = & \prob{(X_1(t) < (1-\delta)tp)  \vee \cdots \vee (X_b(t) < (1-\delta)tp)} \\
 & \leq & \sum_{i=1}^b \prob{ X_i (t) < (1- \delta)tp} \\
 & < & b \exp\left(-tp\frac{\delta^2}{2}\right)
\end{eqnarray*}
where we have applied Chernoff's inequality. Now we choose time $\tau = \frac{8}{p}(d + \log b)$,  and $\delta = 1/2$. Notice that $\tau \geq 2d/p$ and so $d \leq \tau p/2$ . Therefore
\begin{eqnarray*}
\prob{Y(\tau) < d} & \leq & \prob{Y(\tau) < \frac{1}{2} \tau p} \\
& <  & b \exp\left(-\frac{\tau p}{8}\right) \\
& = &  b \exp\left(-d -\log b \right) \\
& = & \exp(-d) \\
& \leq & \exp(-1)
\end{eqnarray*}
So the probability that we have achieved the desired state by time $\tau$ is at least $1 - e^{-1}$.
Applying the lemma, we conclude that the expected time to completion is less than 
\[
 \frac{8(d + \log b)}{p(1 - e^{-1})}
\]
\hfill
$\Box$
 
The upper bound for the star is also an upper bound for the original network. Since the eccentricity of any node in a network is less than the diameter $D$ of the network (the length of the greatest distance between nodes of the network), and the number of leaf nodes $b$ is less than $n$, we have a general upper bound on the propagation time for networks of $O((D + \log n)/p)$. We interpret the diameter $D$ as the time associated with the \emph{depth} of the network, and the factor $\log n$ with the \emph{breadth}. The bound is the maximum of these two factors. 
 
 \section{Results for various networks}
 
We can use these bounds to derive asymptotic results for a range of network topologies. A \emph{random} graph on $n$ nodes is created by assigning an edge between nodes with a given probability. The diameter of such graphs grows as $\log n$. Applying our bound then gives a propagation time of $\Theta((\log n)/p)$. 
 
 For scale-free networks with degree distribution $p(k) \propto k^{-\lambda}$ there are three cases~\cite{cohen:scalefree}. For $\lambda > 3$, the diameter grows as $\log n$ and so again the propagation time is $\Theta((\log n)/p)$. For $2 < \lambda < 3$, the diameter grows much more slowly, as $\log \log n$. In this case the propagation time is between $\log \log n$ and $\log n$. The third case is $\lambda = 3$ for which the diameter grows as $\log n / \log \log n$, which again acts as a lower bound on the propagation time.

Hierarchies in organisational structures may be modelled by tree networks. A complete binary tree has depth $\log n$, and so the propagation time, starting at the root, node is $\Theta(\log n)$. A lattice structure is commonly used in artificial life models (such as cellular automata). Each individual is connected to the neighbours to the north, east, south and west. The diameter of such a network (and therefore the propagation time) is $\Theta(\sqrt{n})$. This is considerably slower than for random and small-world networks. It is known that small-world networks can be constructed from lattices by introducing a small number of random ``short-cuts'' between nodes~\cite{watts:2004}. We see that by doing this, we dramatically reduce the propagation time.
  
If the diameter grows logarithmically in the number of nodes or faster, then it determines the propagation time (that is, it dominates the ``breadth'' factor given by the number of branches in the spanning tree). In this case, the propagation time is also inversely proportional to $p$. However, if it grows slower than logarithmically, we do not get so much information from our bounds. For example, both the complete graph and the hub have constant diameters: our bounds cannot distinguish these cases.
  
 \begin{figure}
\begin{center}
\includegraphics[scale=0.5]{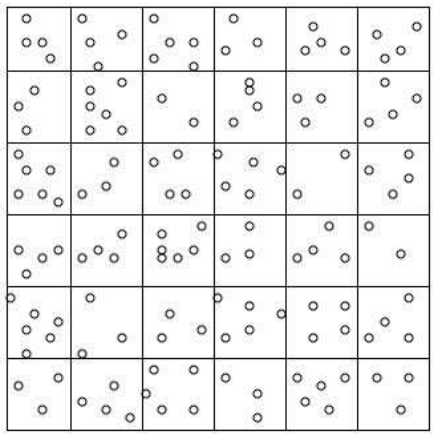}
\end{center}
\caption{Nodes are distributed randomly in the unit square and are connected if they are within a distance of $r$ from each other. We tile the square so that the nodes in each tile form a complete graph.}
\label{fig4}
\end{figure}

This situation occurs in our final example, in which nodes are spatially distributed. Imagine a square with unit length sides. Nodes are distributed randomly in the square and an edge is drawn between nodes that are less than a distance $r$ apart. For example, this could model a random distribution of processors in a computational lawn that have a limited transmission range~\cite{abelson:2000}. The furthest apart two nodes can be geometrically is $\sqrt{2}$, so the diameter of the network, as $n$ increases, approaches $\sqrt{2}/r$. A lower bound on the propagation time is therefore $c/r$, for some constant $c$, which does not depend on $n$. We will show that the propagation time is also bounded above by a constant and is inversely proportional to $r$. To do this, we divide the square up into disjoint \emph{tiles} with side length $r/\sqrt{2}$ (see figure~\ref{fig4}). The diagonal of each tile is $r$, so all the nodes in a tile are connected to all the others. The nodes of one tile in isolation form a complete graph, for which the propagation time is constant. Now consider two neighbouring tiles, with a common edge. Place a third tile so that it covers half of each of these. The nodes in the third tile again form a complete graph, with constant propagation time. This means that the expected time for propagation from one of the original tiles to its neighbour is a constant. The situation, therefore, reduces to the constant time spread of information from tile to tile. Since there are $\sqrt{2}/r$ tiles along each side of the square, the number of tiles that have to be traversed on a path between the corners is proportional to $1/r$. Since each move takes place in constant time, the result follows.

\subsection*{Acknowledgements}
This work was funded by the European Union FP6 project ÒDigital Business EcosystemsÓ.  Dr. Leslie Goldberg (University of Warwick) assisted with the general upper bound.

\bibliographystyle{plain}
\bibliography{propagation}

\end{document}